\documentclass[12pt, verbatim]{amsart}
\usepackage{amscd}
\usepackage{amsthm}
\usepackage{amsmath}
\usepackage{amssymb}
\usepackage{amsbsy}
\usepackage{amsfonts}
\usepackage{euscript}
\usepackage[francais]{babel}
\theoremstyle{definition}
\newtheorem{Thm}{Th\'eor\`eme}
\newtheorem{Prop}{Proposition}
\newtheorem{Lem}{Lemme}
\newtheorem{Cor}{Corollaire}
\newtheorem*{Rq}{Remarque}
\newtheorem*{Rqs}{Remarques}

\newtheorem{Def}{D\'efinition}
\newtheorem{Conj}{Conjecture}
\newtheorem*{Not}{Notation}
\newtheorem*{Nots}{Notations}
\newtheorem*{Conv}{Convention}

\newtheorem{Ex}{Exemple}

\newtheorem{Pete}{Propri\'et\'e}
\newtheorem*{CX_1}{La construction de $X_1$}
\newtheorem*{Hypo}{L'hypoth\`ese de r\'ecurrence}
\newtheorem*{CX_n+1}{La construction de $X_{n+1}$}

\newtheorem*{Ch_n}{La construction de $h_n$}
\newtheorem{Cas}{Cas}

\begin{document}
\title{Un $3-$polyGEM de cohomologie modulo 2 nilpotente}
\author{DongHua JIANG}
\address{LAGA, Institut Galil\'ee, Universit\'e Paris Nord, 93430 Villetaneuse, France}
\email{donghua.jiang@polytechnique.org}
\date{\today}
\maketitle

\begin{abstract}

En 1983, C. McGibbon et J. Neisendorfer \cite{McN} ont
d\'emontr\'e une conjecture de J.-P. Serre \cite{Serre} montrant
qu'un complexe fini $1$-connexe de cohomologie modulo $2$
non-triviale a de la $2$-torsion dans une infinit\'e de groupes
d'homotopie. En 1985, une autre preuve a \'et\'e donn\'ee
par J. Lannes et L. Schwartz \cite{LS86In}. Ce r\'esultat
a sugg\'er\'e une conjecture plus g\'en\'erale: si la cohomologie
modulo 2 r\'eduite d'un polyGEM 1-connexe quelconque est de type fini
et si elle n'est pas r\'eduite \`a $\{ 0 \}$, alors elle contient
au moins un \'el\'ement non nilpotent. Les r\'esultats de Y.
F\'elix, S. Halperin, J.-M. Lemaire et J.-C. Thomas \cite{FHLT} en
1987, de J. Lannes et L. Schwartz \cite{LS89Is} en 1988, et de J.
Grodal \cite{Grodal} en 1996 la soutenaient.

Dans cet article, on construit un contre-exemple.
\end{abstract}

\section{Introduction}

Par convention dans cet article, la cohomologie modulo 2, i.e.,
sur le corps ${\Bbb F}_2$ sera not\'ee ${\rm H}^\ast X$, et la
cohomologie r\'eduite modulo 2 sera not\'ee
$\tilde{\rm H}^\ast X$. J.-P. Serre a d\'emontr\'e en 1953 le
th\'eor\`eme suivant:

\begin{Thm}
(Serre, \cite{Serre}) Soit $X$ un espace simplement connexe de
type fini en 2. On suppose que:
\begin{itemize}
\item la cohomologie $\tilde{\rm H}^\ast X$ est non triviale;
\item les groupes ${\rm H}^n X$ sont nuls pour tout $n$ assez grand.
\end{itemize}
Alors, pour une infinit\'e d'entiers $n$, la multiplication par $2$ du groupe $\pi_n X$ dans lui-m\^eme n'est pas un isomorphisme. $\hfill{\Box}$
\end{Thm}

\begin{Conv}
On dira qu'un espace $X$ est {\bf de type fini en 2}, si sa
cohomologie (modulo 2) est de dimension finie en chaque degr\'e.
\end{Conv}

Tous les espaces que l'on consid\`erera auront cette propri\'et\'e.

Serre conjecturait que sous les hypoth\`eses du th\'eor\`eme, il
existe une infinit\'e d'entiers $n$ tels que le groupe $\pi_n X$
contient un \'el\'ement non trivial d'ordre 2.

Cette conjecture a \'et\'e d\'emontr\'ee par C. McGibbon et J.
Neisendorfer \cite{McN} en 1983, puis une autre preuve a \'et\'e
donn\'ee par J. Lannes et L. Schwartz \cite{LS86In} en 1985.
Depuis, plusieurs g\'en\'eralisations du th\'eor\`eme de Serre ont
\'et\'e donn\'ees par Y. F\'elix, S. Halperin, J.-M. Lemaire et
J.-C. Thomas \cite{FHLT} en 1987, et par J. Lannes et L. Schwartz 
\cite{LS89Is} en 1988. Avant de donner ces \'enonc\'es, rappelons 
quelques d\'efinitions.

\begin{Def}
On dira qu'un espace poss\`ede une {\bf tour de Postnikov finie en
$\bf p = 2$} si la multiplication par 2 du groupe $\pi_n X$ dans
lui-m\^eme est un isomorphisme pour tout $n$ assez grand.
\end{Def}

\begin{Def}\label{polyGEM}
L'ensemble des {\bf polyGEMs} ou {\bf syst\`emes de Postnikov}
(resp. {\bf polyGEMs stables}) est d\'efini r\'ecursivement comme
suit:

\begin{itemize}
\item les $1-$polyGEMs sont les espaces
$$\prod_{1 \leq i \leq i_m, m \in {\Bbb N}} {\rm K}({\rm k}_{i,m}, m),$$
o\`u ${\rm k}_{i,m} \in \{ {\Bbb Z}, {\Bbb Z}/2^l, l = 1, 2,
\cdots, \infty \}$ et $i_m \in {\Bbb N}$ pour tout $m \geq 1$, les
$1-$polyGEMs sont tous stables;

\item les $n-$polyGEMs (resp. les
$n-$polyGEMs stables, qui sont des H-espaces) sont
obtenus comme fibres homotopiques d'une application $f: E \to B$
telle que $E$ est un $(n-1)-$polyGEM (resp. $(n-1)-$polyGEM
stable) et que $B$ est un $1-$polyGEM (dans le cas stable, $f$ est un
morphisme de H-espaces).
\end{itemize}
\end{Def}

\begin{Rqs}
1. Dans la suite, on dira $''$polyGEM$''$ pour $''k-$polyGEM$''$ s'il
n'y a pas lieu de sp\'ecifier l'entier $k$.

2. La d\'efinition ci-dessus n'est pas la plus g\'en\'erale (voir
\cite{Farjoun}), on pourrait d\'efinir un $1-$polyGEM par la
condition que c'est un produit (infini) d'espaces
d'Eilenberg-MacLane ${\rm K}(G_n,n)$, $n = 1, 2, \cdots$. Puis
proc\'eder it\'erativement comme plus haut, modulo quelques
hypoth\`eses tout polyGEM est homotopiquement \'equivalent \`a
$2$-compl\'etion pr\`es \`a un de ceux dans la d\'efinition
\ref{polyGEM}.
\end{Rqs}

On a:

\begin{Thm}\label{Lannes}
(Lannes et Schwartz, \cite{LS89Is}, 0.1) Soit $X$ un espace
simplement connexe de type fini en 2. On suppose que:
\begin{itemize}
\item la cohomologie $\tilde{\rm H}^\ast X$ est non triviale;

\item la cohomologie $\tilde{\rm H}^\ast X$ est nilpotente
({\it i.e., tout \'el\'ement de $\tilde{\rm H}^\ast X$ est nilpotent}).
\end{itemize}
Alors, pour une infinit\'e d'entiers $n$, la multiplication par
2 du groupe $\pi_n X$ dans lui-m\^eme n'est pas un isomorphisme.
Autrement dit, $X$ n'a pas de tour de Postnikov finie en $p=2$. $\hfill{\Box}$
\end{Thm}

\begin{Rq}
Si on suppose, dans la d\'efinition des polyGEMs, que les $i_m$
sont \'egaux \`a z\'ero d\`es que $m$ est assez grand, on obtient
des espaces que l'on appellera les {\bf polyGEMs finis} ({\bf en $p =
2$}), ce sont des espaces dont la tour de Postnikov est finie en 2.
Alors le th\'eor\`eme \ref{Lannes} est \'equivalent \`a dire que
la cohomologie r\'eduite des polyGEMs finis simplement connexes est soit 
triviale, soit non nilpotente.
\end{Rq}

On a aussi:

\begin{Thm}\label{Felix}
(F\'elix, Halperin, Lemaire et Thomas, \cite{FHLT}, 5.1) Soit $X$
un espace simplement connexe de type fini en 2. On suppose que:
\begin{itemize}
\item la cohomologie $\tilde{\rm H}^\ast X$ est non triviale;
\item ${\rm cat} X$ est fini ({\it i.e., $X$ poss\`ede un recouvrement fini par des espaces contractiles et on note par ${\rm cat} X$ le nombre minimum d'espaces pour un tel recouvrement}).
\end{itemize}
Alors, $X$ n'est pas un polyGEM. $\hfill{\Box}$
\end{Thm}

De plus, d'apr\`es J. Grodal:

\begin{Thm}\label{Grodal}
(Grodal, \cite{Grodal}, 6.5) Soit $X$ un polyGEM simplement
connexe de type fini en 2. Si la cohomologie
$\tilde{\rm H}^\ast X$ est non triviale, alors elle
n'est pas localement finie en tant que module sur l'alg\`ebre de
Steenrod. $\hfill{\Box}$
\end{Thm}

Voici quelques pr\'ecisions sur cet \'enonc\'e.

\begin{itemize}
\item Un ${\mathcal A}_2-$module est {\bf localement fini} si et
seulement si le sous-module sur l'alg\`ebre de Steenrod engendr\'e
par tout \'el\'ement est fini.

\item Si un ${\mathcal A}_2-$module n'est pas nilpotent, il n'est
pas localement fini non plus. Au contraire, un ${\mathcal
A}_2-$module qui n'est pas localement fini peut \^etre nilpotent,
voir par exemple la ${\mathcal A}_2-$alg\`ebre ${\mathcal E}(n)$
d\'efinie dans la prochaine section.
\end{itemize}

Dans son article \cite{Grodal}, Grodal propose:

\begin{Conj}\label{GrodalConj}
(Grodal, \cite{Grodal}, 6.1)
La cohomologie r\'eduite des polyGEMs simplement connexes de type 
fini en 2 est soit triviale, soit non nilpotente.
\end{Conj}

Cette conjecture semble raisonnable, \'etant donn\'es les r\'esultats
pr\'ec\'edents. Mais si cet \'enonc\'e est faux, en particulier
pour les polyGEMs stables, alors il existe un polyGEM stable $X$
dont la cohomologie $\tilde{\rm H}^\ast X$ est non triviale et
nilpotente. Le th\'eor\`eme \ref{Lannes} dit que pour une
infinit\'e d'entiers $n$, la multiplication par 2 du groupe $\pi_n
X$ dans lui-m\^eme n'est pas un isomorphisme, et le th\'eor\`eme
\ref{Felix} dit que ${\rm cat} X = + \infty$, enfin celui de
\ref{Grodal} affirme que la cohomologie n'est pas localement
finie.

L'objectif de cet article est de donner un contre-exemple \`a
cette conjecture. Etant donn\'e un entier $l \geq 2$, on va construire un
$3-$polyGEM stable $(2l-1)-$connexe $X$ dont la cohomologie
$\tilde{\rm H}^\ast X$ est non triviale et nilpotente.

Le plan de cet article est la suivante: la section 2 contient des
pr\'eliminaires, i.e., des travaux de J. Milgram et de L. Smith.
On construit les espaces $X_n$ dans la section 3 dont la $''$limite$''$
$X_\infty$ est le contre-exemple cherch\'e. En revanche dans la
derni\`ere section, on montre que la conjecture \ref{GrodalConj}
est vraie pour les $2-$polyGEMs stables simplement connexes.

\section{Pr\'eliminaires}

Cette section rappelle des travaux de J. Milgram et de L. Smith. Le
travail de Milgram donne la structure de la cohomologie d'un
certain $2-$polyGEM en tant qu'alg\`ebre de Hopf, avec des
informations sur la structure de module instable. Celui de Smith
donne essentiellement la structure d'alg\`ebre, \`a filtration
pr\`es, pour tous les polyGEMs stables.

Rappelons que:
\begin{itemize}
\item {\bf le module $\bf {\rm F}(n)$} est le module instable libre sur
l'alg\`ebre de Steenrod engendr\'e par un g\'en\'erateur $\iota_n$
de degr\'e $n$, une base sur le corps ${\Bbb F}_2$, \`a suspension
$n$-i\`eme pr\`es, est donn\'ee par l'ensemble des op\'erations
de Steenrod d'exc\`es inf\'erieur ou \'egal \`a $n$;

\item {\bf l'op\'eration $\bf Sq_0$} est d\'efinie sur un module
instable $M$ par $Sq_0 x = Sq^{|x|} x$ pour tout $x \in M$. En
particulier, l'ensemble $Sq_0 M$ est un sous-module instable de
$M$:
$$Sq_0 M = \{ Sq_0 x\ |\ x \in M \}$$
et $Sq_0 {\rm F}(n)$ s'identifie au ${\Bbb F}_2-$espace vectoriel
gradu\'e engendr\'e par les $Sq^I \iota_n$, ${\rm ex}(I) = n$;

\item et {\bf le foncteur $\bf U$ de Steenrod-Epstein}
\cite{Steenrod} est le foncteur qui associe \`a tout module
instable $M$ l'alg\`ebre (instable) enveloppante
$${\bf U}(M) = S^\ast(M) / (Sq_0 x - x^2, x \in M),$$
o\`u, $x^2$ d\'esigne le carr\'e dans l'alg\`ebre sym\'etrique.
C'est l'adjoint \`a gauche du foncteur {\it oubli} de la
cat\'egorie des alg\`ebres instables vers celle des modules
instables, ${\bf U}(M)$ est une alg\`ebre de Hopf primitivement
engendr\'ee \cite{MM};

\item la cohomologie ${\rm H}^{*} ({\rm K}({\Bbb F}_2, n))$ est
isomorphe \`a ${\bf U}({\rm F}(n))$ en tant qu'alg\`ebre de Hopf.
Comme alg\`ebre polyn\^omiale sur le corps ${\Bbb F}_2$ elle est
engendr\'ee par les classes $Sq^I \iota_n$ avec $I$ suite
admissible d'exc\`es $< n$.
\end{itemize}

\begin{Def}
(Milgram, \cite{Milgram}, 1.1.1) Soit $E_n$ la fibre homotopique
de l'application ${\rm K}({\Bbb F}_2, n) \to {\rm K}({\Bbb F}_2,
2n)$ qui correspond au cup-carr\'e de la classe (primitive)
$\iota_n \in {\rm H}^n ({\rm K}({\Bbb F}_2, n))$. L'espace $E_n$
est un $2-$polyGEM stable et donc un H-espace, ${\rm H}^\ast E_n$
est une alg\`ebre de Hopf.
\end{Def}

\begin{Rq}
Comme remarqu\'e par F. Cohen \cite{Cohen}, ces espaces \'etaient d'abord \'etudi\'es par L. Kristensen \cite{Kristensen}, et puis par E.H. Brown et F.P. Peterson \cite{BP} lorsqu'ils calculaient les cup-produits en bas degr\'e qui permettent de construire des op\'erations cohomologiques instables secondaires d\'etectant le produit de Whitehead (sur les sph\`eres de dimension diff\'erente de $2^k-1$).
\end{Rq}

\begin{Conv}
On appellera l'espace $E_n$ un {\bf espace de Milgram} dans la suite.
\end{Conv}

\begin{Nots}
1. $\bar{\iota}_n$ d\'esignera un \'el\'ement de degr\'e $n$ tel
que, pour toute suite admissible $I$ d'exc\`es $< n$ les
\'el\'ements $Sq^I \bar{\iota}_n$ sont lin\'eairement
ind\'ependants, et que $\bar{\iota}_n^2 = 0$.
Le ${\mathcal A}_2-$module engendr\'e est donc le quotient
${\rm F}(n) / Sq_0 {\rm F}(n)$.

2. $\lambda$ d\'esignera g\'en\'eriquement un ${\mathcal
A}_2-$g\'en\'erateur primitif (d'une ${\mathcal A}_2-$alg\`ebre de
Hopf); $\tau$ d\'esignera g\'en\'eriquement un ${\mathcal
A}_2-$g\'en\'erateur non primitif (d'une ${\mathcal
A}_2-$alg\`ebre de Hopf).
\end{Nots}

Par d\'efinition l'alg\`ebre instable ${\mathcal E}(n)$ est
${\bf U}({\rm F}(n) / Sq_0 {\rm F}(n))$.
En tant qu'alg\`ebre c'est  l'alg\`ebre
ext\'erieure engendr\'ee par les $Sq^I \bar{\iota}_n$, $I$ suite
admissible d'exc\`es strictement inf\'erieur \`a $n$.

Soit ${\mathcal I}_{2n-1}$ le sous-module de ${\rm F}(2n-1)$ engendr\'e
par les $Sq^K {\iota}_{2n-1}$ o\`u $K = (k_1, \cdots, k_r)$ est
une suite admissible (d'exc\`es strictement inf\'erieur \`a $2n$)
telle qu'il existe un $k_i$ impair, $1 \leq i \leq r$. Par
d\'efinition ${\mathcal P}(n) = {\bf U}({\mathcal I}_{2n-1})$.
Notons que $\iota_{2n-1}$ n'est pas dans ${\mathcal I}_{2n-1}$.

Le th\'eor\`eme de Milgram dit que la cohomologie (modulo 2)
de $E_n$ est isomorphe \`a ${\mathcal E}(n)
\otimes {\mathcal P}(n)$ en tant qu'alg\`ebre. Il pr\'ecise la
diagonale.

\begin{Thm}\label{Milgram}
(Milgram, \cite{Milgram}, 1.2.1, 1.3.3)
On a un isomorphisme d'alg\`ebres
$${\rm H}^\ast E_n \cong {\mathcal E}(n) \otimes {\mathcal P}(n).$$
\begin{itemize}
\item pour $n$ pair, en tant que ${\mathcal A}_2-$alg\`ebre,
${\rm H}^\ast E_n$ poss\`ede les g\'en\'erateurs
$\bar{\iota}_n$, $\tau_1$, $\lambda_2$, $\cdots$, $\lambda_k$ dont
$\bar{\iota}_n$, $\lambda_2$, $\cdots$, $\lambda_k$ sont primitifs
et
$$\Delta(\tau_1) = \tau_1 \otimes 1 + \bar{\iota}_n \otimes \bar{\iota}_n + 1 \otimes \tau_1;$$
\item pour $n$ impair, en tant que ${\mathcal A}_2-$alg\`ebre,
${\rm H}^\ast E_n$ poss\`ede les g\'en\'erateurs
$\bar{\iota}_n$, $\lambda_1$, $\tau_2$, $\cdots$, $\tau_k$ dont
$\bar{\iota}_n$, $\lambda_1$ sont primitifs et pour $i = 2,
\cdots, k$,
$$\Delta(\tau_i) = \tau_i \otimes 1 + Sq^{2^{i-1}-1}
\bar{\iota}_n \otimes Sq^{2^{i-1}-1} \bar{\iota}_n + 1 \otimes
\tau_i.$$
\end{itemize}
O\`u, $k$ est l'unique entier tel que $2^{k-1} \leq n < 2^k$, et
$|\tau_i| = |\lambda_i| = 2(n+2^{i-1}-1)$, $\forall$ $i = 1,
\cdots, k$. $\hfill{\Box}$
\end{Thm}

\begin{Rqs}
1. Les $\bar{\iota}_n$, $\tau_i^2$ et $\lambda_j$, $i, j = 1, \cdots,
k$, sont primitifs.

2. Dans le m\^eme papier, Milgram donne la structure de module sur
l'alg\`ebre de Steenrod pour un autre syst\`eme de
g\'en\'erateurs, \'equivalent \`a celui dans l'\'enonc\'e. Et il
pr\'ecise le lien entre ces deux syst\`emes de g\'en\'erateurs.
\end{Rqs}

On rappelle maintenant les travaux de Smith \cite{Smith}. En fait le
r\'esultat de Milgram se d\'eduit pour partie au moins de ceux-ci.
Plus pr\'ecis\'ement, on va calculer, \`a l'aide de la suite
spectrale d'Eilenberg-Moore, la cohomologie de la fibre $F$
associ\'ee \`a une fibration $\pi: E \to B$, o\`u $B$ est un
1-polyGEM 1-connexe et $E$ est un H-espace. On a un carr\'e
cart\'esien:
\begin{equation}\label{carre1}
\begin{CD}
F @>>> \ast \\
@VVV   @VVV \\
E @>{\pi}>> B
\end{CD}
\end{equation}
Comme dans \cite{Smith}, on suppose que $\pi$ est un morphisme de
H-espaces pour les structures de H-espace de $E$ et de $B$, et que
${\rm H}^\ast E$ est une alg\`ebre de Hopf cocommutative.
On remarque que par ces hypoth\`eses, ${\rm ker}(\pi^*)$ est un
id\'eal de Hopf dans ${\rm H}^\ast B$.

\begin{Not}
$\{ {\rm E}_r, d_r \}$ d\'esignera la suite spectrale d'Eilenberg-Moore associ\'ee au carr\'e cart\'esien (\ref{carre1}) v\'erifiant les hypoth\`eses dans le paragraphe pr\'ec\'edent.
\end{Not}

\begin{Thm}\label{Thm2}
(Smith, \cite{Smith}, 1.5) Soient $\Gamma$, $A$ deux alg\`ebres de
Hopf cocommutatives et soit $\varphi: \Gamma \to A$ un
morphisme d'alg\`ebres de Hopf. Soit
$\Lambda = {\rm sub}$-${\rm ker}(\varphi)$,
la sous-alg\`ebre de Hopf de $\Gamma$ qui engendre
${\rm ker}(\varphi)$. Alors il y a un isomorphisme d'alg\`ebres de Hopf:
$${\rm Tor}_\Gamma (A, {\Bbb F}_2) \cong A//\varphi \otimes {\rm Tor}_\Lambda ({\Bbb F}_2, {\Bbb F}_2).$$ $\hfill{\Box}$
\end{Thm}

\begin{Cor}\label{Cor2}
(cf. \cite{Smith}, 2.1)
On a un isomorphisme d'alg\`ebres:
$${\rm E}_2 \cong {\rm H}^\ast E // {\rm im}(\pi^\ast) \otimes {\rm Tor}_{{\rm sub-ker}(\pi^\ast)}^\ast ({\Bbb F}_2, {\Bbb F}_2).$$
\end{Cor}

\begin{proof}
Par hypoth\`ese, $\pi^\ast$ est un morphisme d'alg\`ebres de Hopf et ${\rm H}^\ast E$ est une alg\`ebre de Hopf cocommutative.
D'autre part, le calcul de Serre (voir \cite{Serre}) dit que ${\rm H}^\ast B$ est une alg\`ebre de Hopf cocommutative. Donc, on peut \'etablir le r\'esultat d'apr\`es le th\'eor\`eme \ref{Thm2}.
\end{proof}

\begin{Prop}\label{Prop3}
(cf. \cite{Smith}, 2.2)
${\rm E}_2 = {\rm E}_\infty$.
\end{Prop}

\begin{proof}
Le calcul de Serre (voir \cite{Serre}) dit que ${\rm H}^\ast B =
{\rm P}[V]$, i.e., une alg\`ebre polyn\^omiale en certain espace
vectoriel gradu\'e $V$. Comme ${\rm sub}$-${\rm ker}(\pi^\ast)$
est une sous-${\mathcal A}_2$-alg\`ebre de Hopf de ${\rm P}[V]$,
le th\'eor\`eme de Borel (voir \cite{MM}, 7.11) sur la structure
des alg\`ebres de Hopf sur le corps ${\Bbb F}_2$ dit que c'est
aussi une alg\`ebre polyn\^omiale, i.e., ${\rm sub}$-${\rm
ker}(\pi^\ast) = {\rm P}[x_1, \cdots, x_n, \cdots]$. Plus
pr\'ecis\'ement elle s'identifie \`a ${\bf U}(P)$, o\`u $P$ est le
sous-module instable des \'el\'ements primitifs de la sous-alg\`ebre
de Hopf ${\rm sub}$-${\rm ker}(\pi^\ast)$.

A l'aide du complexe de Koszul, on obtient:
$${\rm Tor}_{{\bf U}(P)}^\ast ({\Bbb F}_2, {\Bbb F}_2) = {\rm E}[(P/Sq_0(P))].$$
Soit, comme espace vectoriel gradu\'e, la base $u_i$, $i \in I$,
du module instable $P/Sq_0(P)$ qui correspond \`a la base monomiale des $x_i$ de $P$. En tant qu'alg\`ebre, on a avec un petit abus de notation:
$${\rm E}_2 \cong {\rm H}^\ast E // {\rm im}(\pi^\ast) \otimes {\rm E}[u_1, \cdots, u_n, \cdots].$$
Dans ces formules ${\rm E}[\cdots]$ d\'esigne l'alg\`ebre
ext\'erieure, soit sur l'espace vectoriel gradu\'e $P/Sq_0(P)$ en
degr\'e cohomologique $-1$, soit sur les g\'en\'erateurs $u_i$ en
bidegr\'e $(-1, |x_i|)$.

Donc ${\rm E}_2$ est une alg\`ebre engendr\'ee par ${\rm E}_2^{0,
\ast}$ et ${\rm E}_2^{-1, \ast}$.

Or la diff\'erentielle $d_r$ agit comme une d\'erivation, et pour $p = 0$ ou 1, $r \geq 2$, on a
$$d_r: {\rm E}_2^{-p, \ast} \to {\rm E}_2^{-p+r, \ast} = 0.$$
Donc pour des raisons de degr\'e, $d_r = 0$, $\forall$ $r \geq 2$
et ${\rm E}_2 = {\rm E}_\infty$.
\end{proof}

Etant donn\'ee une fibration $\pi: E \to B$ de fibre
$F$, la cohomologie ${\rm H}^\ast F$ a une filtration
d\'ecroissante par des modules instables $F_s = F_s ({\rm H}^\ast
F)$, $s \leq 0$ et $F_0$ est l'image de ${\rm H}^\ast E$. De plus
le produit envoie $F_s \otimes F_t$ vers $F_{s+t}$, l'objet gradu\'e
associ\'e ${\bf Gr}\, {\rm H}^\ast F$ a donc une structure d'alg\`ebre et de module instable et:

\begin{Cor}\label{Gr}
Etant donn\'ee une fibration $\pi: E \to B$ de fibre $F$, o\`u $B$ est un
1-polyGEM 1-connexe, $\pi$ est un morphisme de H-espaces et ${\rm
H}^\ast E$ est une alg\`ebre de Hopf cocommutative. On a un
isomorphisme de ${\mathcal A}_2-$alg\`ebres gradu\'ees:
$${\bf Gr}\, {\rm H}^\ast F \cong {\rm H}^\ast E // {\rm im}(\pi^\ast) \otimes {\rm E}[u_1, \cdots, u_n, \cdots].$$ $\hfill{\Box}$
\end{Cor}

\begin{Rqs}
1. Dans la formule ci-dessus le degr\'e de $u_i$ consid\'er\'ee comme
classe dans ${\bf Gr}\ {\rm H}^\ast F $ est $|x_i|-1$.

2. Ce corollaire r\'esulte du corollaire \ref{Cor2} et de la proposition
\ref{Prop3}. En effet, on utilise l'id\'ee de Smith (voir
\cite{Smith}, $\S 2$).

3. Un moment de r\'eflexion montre que ceci
restitue le r\'esultat de Milgram pour la partie multiplicative.
\end{Rqs}

Pr\'ecisons un peu le r\'esultat, en tant que module instable on a:
$$\Sigma^{-s}F_s/F_{s+1} \cong {\rm E}_\infty^{s, \ast}.$$
Donc dans notre cas:
$$\Sigma^{-s}F_s/F_{s+1} \cong {\rm E}_2^{s, \ast},$$
et
$$\Sigma^{-s}F_s/F_{s+1} \cong {\rm H}^\ast E // {\rm im}(\pi^\ast)
\otimes E_i,$$
o\`u $E_i$ est la $i$-\`eme puissance ext\'erieure de
$\Sigma^{-1}P/Sq_0(P)$, donc est engendr\'ee par les mon\^omes
$u_{a_1} \cdots u_{a_i}$.

\begin{Ex}
Pour mieux comprendre que le morphisme dans le corollaire \ref{Gr}
n'est pas un isomorphisme d'alg\`ebres en g\'en\'eral, on \'etudie
la fibration
$$\begin{CD}
F @>>> E @>{\pi}>> B
\end{CD}$$
avec $B = {\rm K}({\Bbb F}_2, 2)$, $E = {\rm P}B \cong \ast$ et $F = \Omega B \cong {\rm B}{\Bbb F}_2$. Comme $\tilde{\rm H}^\ast E = \{ 0 \}$ et ${\rm ker}(\pi^\ast) = {\rm H}^\ast B$ est une alg\`ebre de Hopf, on a donc
$$\begin{array}{rcl}
\text{sub-ker}(\pi^\ast)& =& {\rm ker}(\pi^\ast) \\
& =& {\rm H}^\ast B \\
& =& {\rm P}[Sq^{2^{i-1}} \cdots Sq^2 Sq^1 \iota_2, i = 0, 1, \cdots]
\end{array}$$
et
$$\begin{array}{l}
{\rm H}^\ast E // {\rm im}(\pi^\ast) \otimes {\rm E}[u_1, \cdots, u_n, \cdots] \cong \\
\qquad \qquad \qquad \qquad {\rm E}[\Sigma^{-1}\{ Sq^{2^{i-1}} \cdots Sq^2 Sq^1 \iota_2, i = 0, 1, \cdots \} ]
\end{array}$$
est une alg\`ebre ext\'erieure (gradu\'ee). Or ${\rm H}^\ast F = {\Bbb F}_2 [ \iota_1 ]$ est une alg\`ebre polyn\^omiale, le morphisme en question n'est donc pas un isomorphisme d'alg\`ebres.
\end{Ex}

Avant de terminer cette section on \'enonce un lemme qui sera utile dans
la suite. Il donne des informations sur les g\'en\'erateurs polyn\^omiaux de la
cohomologie de la fibre $F$ d'une application de H-espaces $\pi: E
\to B$, o\`u $B$ est un $1-$polyGEM $1-$connexe. La
cohomologie de $E$ est suppos\'ee cocommutative.

\begin{Lem}\label{Lemme1}
Supposons que le noyau de $\pi^\ast$, en tant qu'alg\`ebre de
Hopf, soit engendr\'e par des cup-carr\'es. Alors pour tout
\'el\'ement non-nilpotent $x$ de ${\rm H}^\ast F$, il existe un
entier positif $k$ tel que $Sq_0^k (x)$ soit dans l'image de ${\rm
H}^\ast E$.
\end{Lem}

\begin{proof}
Soit, comme plus haut, $P$ l'ensemble des \'el\'ements primitifs
dans le noyau de $\pi^\ast$. Alors les r\'esultats de Smith, en
particulier le corollaire \ref{Gr} (rappellons que $B$ est un 
$1-$polyGEM $1-$connexe), disent qu'un gradu\'e comme
module instable de ${\rm H}^\ast F$ est isomorphe \`a:
$${\rm H}^\ast E //{\rm im} (\pi^\ast) \otimes {\rm E}[\Sigma^{-1}P/Sq_0(P)].$$
Par hypoth\`ese, $P \cong Sq_0M$ pour certain module instable $M$. 
Par cons\'equent, le module $\Sigma^{-1}P/Sq_0(P)$ est une suspension.

Soit $x$ un \'el\'ement non nilpotent de ${\rm H}^\ast F$. Tous
les \'el\'ements en degr\'e strictement positif de ${\rm
E}[\Sigma^{-1}P/Sq_0(P)]$ sont nilpotents.  Supposons que la
classe $\bar x \not = 0$ de $x$ dans le gradu\'e soit de la forme
$\sum_\ell \, p_\ell \otimes y_\ell \not = 0$, $y_\ell$
appartenant \`a la $i$-\`eme puissance ext\'erieure pour un entier
$i > 0$. La classe $x$ est donc dans le terme $F_{-i}$ de la
filtration de ${\rm H}^\ast F$ mais pas dans $F_{-i+1}$. Il existe
$d > 0$ tel que $Sq_0^d(y_\ell) = 0$ pour tout $\ell$. Donc
$$Sq_0^d(\bar{x}) = \sum_\ell Sq_0^d(p_\ell) \otimes Sq_0^d(y_\ell) = 0,$$
i.e., la classe de $Sq_0^d(x)$ est au moins dans le terme $F_{-i+1}$ de la
filtration. Par it\'eration on obtient que $Sq_0^k(x) \in {\rm
H}^\ast E // {\rm im}(\pi^\ast)$ pour certain entier $k > 0$,
l'affirmation suit.
\end{proof}

\section{La construction du contre-exemple}

Dans cette section, on va \'etablir le r\'esultat principal de
l'article:

\begin{Thm}\label{3-polyGEM}
Etant donn\'e un entier $l \geq 2$, il existe un 3-polyGEM stable $(2l-1)-$connexe dont la cohomologie r\'eduite modulo 2 est non triviale et nilpotente.
\end{Thm}

Pour \'etablir ce r\'esultat, on va proc\'eder en plusieurs
\'etapes. Plus pr\'ecis\'ement, on donne d'abord la
construction des espaces $X_n$ dont la $''$limite$''$ $X_\infty$ est
l'espace cherch\'e. Ceci est fait dans la premi\`ere sous-section.
Dans la seconde on \'etablit les propri\'et\'es de cette
construction. Dans la derni\`ere on
montre que l'espace $X_\infty$ est un bon candidat qui permet
d'\'etablir le th\'eor\`eme \ref{3-polyGEM}.

\subsection{La construction des espaces $X_n$}

On va donc construire un $3-$polyGEM stable dont la cohomologie
r\'eduite modulo 2 est nilpotente. Pour ce faire, on construit des
$2-$polyGEMs stables finis $P_n$ (qui sont des produits d'espaces
$E_i$ de Milgram) et $X_n$ est la fibre homotopique d'une
application de $P_n$ dans un 1-polyGEM $G_n$:
$$\begin{CD}
X_n @>{i_n}>> P_n @>{f_n}>> G_n
\end{CD}$$
telle que le noyau (en tant qu'alg\`ebre de Hopf) soit
engendr\'e par des cup-carr\'es. Par construction on aura
des applications $p_{n+1}: X_{n+1} \to X_{n}$, telles que
$$p_{n+1}^\ast: {\rm H}^i X_n \to {\rm H}^i X_{n+1}$$
soit un isomorphisme pour tout entier $i \leq 4n-2$.
Ces applications induiront des
applications $q_n: X_{\infty} \to X_n$, telles que
$$q_n^\ast: {\rm H}^i X_n \to {\rm H}^i X_\infty$$
soit un isomorphisme pour tout entier $i \leq 4n-2$, de plus les images
dans ${\rm H}^\ast X_\infty$ par $q_n^\ast$ de tous les
g\'en\'erateurs (en tant qu'alg\`ebre) de $\tilde{\rm H}^\ast X_n$
seront nilpotentes.

L'espace cherch\'e est obtenu pour $n = \infty$, la condition
cherch\'ee (la nilpotence de tout \'el\'ement en degr\'e
strictement positif) sera cons\'equence de calculs dans la
suite spectrale d'Eilenberg-Moore.
Les espaces $E_i$ de Milgram sont les briques \'el\'ementaires de
cette construction. On va proc\'eder par r\'ecurrence sur $n$ pour
construire l'espace $X_n$, l'id\'ee est de rendre
nilpotents les g\'en\'erateurs (en tant qu'alg\`ebre)
par r\'ecurrence sur le degr\'e.
La proc\'edure standard pour ce faire produit une tour
(et donc un syst\`eme) de Postnikov infinie, avec les espaces
$E_i$ on s'affranchit de cette contrainte.

\begin{CX_1}
Pour $n = 1$, la construction est la suivante: soit $l \geq 2$,
d\'efinissons $X_1 = E_{2l}$, $P_1 = E_{2l}$ et $G_1 =
\ast$. On a la fibration triviale:
$$\begin{CD}
X_1 @>>> P_1 @>{f_1}>> G_1 = \ast,
\end{CD}$$
o\`u ${\rm ker}(f_1^\ast) = \{ 0 \}$ et o\`u $X_1$ est
$(2l-1)-$connexe.
\end{CX_1}

Supposons la construction faite pour l'entier $n$.

\begin{Hypo}
On va supposer que:

\begin{enumerate}
\item $p_n^\ast: {\rm H}^\ast X_{n-1} \to {\rm H}^\ast X_n$
est un isomorphisme en degr\'e inf\'erieur ou \'egal
\`a $4n-2$;

\item $P_n$ est le produit de $P_{n-1}$ par des espaces de Milgram,
et les images des g\'en\'erateurs (en tant qu'alg\`ebre)
de $\tilde{\rm H}^\ast P_{n-1}$ en degr\'e
$n$ via $p_n^\ast \circ i_{n-1}^\ast$
dans $\tilde{\rm H}^\ast X_n$ sont nilpotentes;

\item $(f_n^\ast)^{-1} (Sq_0 ({\rm PH}^\ast P_n)) \subset
Sq_0({\rm PH}^\ast G_n)$, ce qui implique que le noyau (en tant
qu'alg\`ebre de Hopf) est engendr\'e par des cup-carr\'es.
\end{enumerate}

\begin{Rqs}
1. L'ensemble de ces conditions (pour tous les entiers $n$) impliquera que les
images par $q_{n-1}^\ast$ dans $\tilde{\rm H}^\ast X_\infty$ de
tous les g\'en\'erateurs en tant qu'alg\`ebre de $\tilde{\rm
H}^\ast X_{n-1}$ sont nilpotentes. Par cons\'equent et \`a cause de
la stabilit\'e tout \'el\'ement de degr\'e strictement positif de ${\rm H}^\ast
X_\infty$ est nilpotent. En fait l'image de $x \in {\rm H}^\ast X_{n-1}$,
$|x| > 0$, est nilpotente dans ${\rm
H}^\ast X_{n+k}$ d\`es que $k$ est assez grand.

2. Dans la condition $(3)$ et dans la suite $f_n^\ast$
d\'esignera aussi l'application restreinte aux primitifs,
le contexte fera qu'il n'y aura pas d'ambiguit\'es.

3. La condition $(1)$ est trivialement r\'ealis\'ee pour $n = 1$.
\end{Rqs}
\end{Hypo}

\begin{CX_n+1}
On suppose $X_n$ construit, on passe \`a $X_{n+1}$.

Si tous les g\'en\'erateurs de $\tilde{\rm H}^\ast X_n$ en degr\'e
$n+1$ provenant de $\tilde{\rm H}^\ast P_n$ sont nilpotents,
on d\'efinit $X_{n+1} = X_n$, $P_{n+1} = P_n$, $G_{n+1} = G_n$,
$f_{n+1} = f_n$.

En particulier, et \'etant donn\'ee la structure de ${\rm H}^\ast E_{2l}$, on a
$$X_1 = X_2 = \cdots = X_{4l-1} = E_{2l}.$$

Dans le cas contraire $P_{n+1}$ sera un produit de la forme:
$$P_n \times E = P_n \times E_{2n+2}^{\times \alpha}$$
pour un certain entier $\alpha$. De m\^eme $G_{n+1}$ sera un produit:
$$G_n \times G = G_n \times K_{2n+2}^{\times \alpha}$$
pour le m\^eme entier $\alpha$, avec $K_{2n+2} = {\rm K}({\Bbb F}_2, 2n+2)$.

Par exemple $X_{4l}$ sera la fibre de l'application de H-espaces:
$$E_{2l} \times E_{8l} \to K_{8l}$$
d\'etermin\'ee comme suit. La cohomologie de $E_{2l}$ a un g\'en\'erateur
polyn\^omial $\tau$ en degr\'e $4l$.
Ce g\'en\'erateur $\tau$ n'est pas primitif mais son carr\'e l'est.
On consid\`ere alors l'application donn\'ee par la
classe $\tau^2+ \bar \iota_{8l}$ qui est primitive. Comme $\bar
\iota_{8l}$ est de carr\'e nul l'image de $\tau$ dans la cohomologie
de la fibre homotopique $X_{4l}$ est de puissance quatri\`eme nulle.

Il faut syst\'ematiser cette construction. L'application
$$f_{n+1}: P_n \times E \to G_n \times G$$
sera de la forme:
$$f_{n+1}(\alpha, \beta) = (f_n(\alpha), h_n(\alpha, \beta)).$$
On aura donc un diagramme commutatif de la forme:
\begin{equation}\label{diag1}
\begin{CD}
X_{n+1} @>>>  P_{n+1} @>{f_{n+1}}>> G_{n+1} \\
@V{p_{n+1}}VV @VVV                  @VVV \\
X_n     @>>>  P_n     @>{f_n}>>     G_n
\end{CD}
\end{equation}
Les deux applications verticales ne portant pas de nom sont les
projections \'evidentes.
\end{CX_n+1}

\begin{Ch_n}
Si le noyau de chaque application $f_n^\ast$ est constitu\'e de
cup-carr\'es, pour garantir la nilpotence il suffit,
gr\^ace au lemme \ref{Lemme1} et \`a la stabilit\'e,
de rendre nilpotente les images des g\'en\'erateurs polyn\^omiaux
de $\tilde{\rm H}^\ast P_n$ dans ${\rm H}^\ast X_n$.
On va donc faire une construction qui rend nilpotente les images
dans ${\rm H}^\ast X_n$, en degr\'e $n+1$, des g\'en\'erateurs
polyn\^omiaux de ${\rm H}^\ast P_n$ (dans le m\^eme degr\'e).
Nous aurons \`a distinguer deux cas selon qu'ils sont
(en tant qu'\'el\'ement de ${\rm H}^\ast P_n$) primitifs ou non.
Rappelons que l'on a distingu\'e dans ${\rm H}^\ast E_m$ des
\'el\'ements de type $\tau$ ou $\lambda$. On notera donc
dans ${\rm H}^\ast P_n$
\begin{itemize}
\item $x_1, \cdots, x_k$ pour les g\'en\'erateurs primitifs;

\item $y_1, \cdots, y_h$ pour ceux qui ne le sont pas,
\end{itemize}
dans ce dernier cas selon Milgram leurs cup-carr\'es le sont.

Posons alors
$$\begin{array}{rcl}
E& =& E_{2n+2}^{\times (k+h)} \\
G& =& K_{2n+2}^{\times (k+h)}
\end{array}$$
et $f_{n+1}$ est d\'efini par la formule suivante:
$$\begin{array}{rrcl}
f_{n+1}: & P_n \times E & \to & G_n \times G \\
& (\alpha, \beta) & \mapsto & (f_n(\alpha), h_n(\alpha, \beta))
\end{array}$$
o\`u il faut d\'efinir $h_n: P_n \times E \to G$, mais comme
$$[P_n \times E, G] \cong \left({\rm H}^{2n+2} (P_n \times E)\right)^{\times (k+h)},$$
il suffit de donner $k+h$ classes de cohomologie, ce seront les
$$x_1^2 + \bar{\iota}_{2n+2,1}, \cdots, x_k^2 + \bar{\iota}_{2n+2,k}$$
et
$$y_1^2 + \bar{\iota}_{2n+2,k+1}, \cdots, y_h^2 + \bar{\iota}_{2n+2,k+h}$$
qui sont toutes primitives, et donc $h_n$ est une application de H-espaces
($\bar{\iota}_{2n+2,i}$ d\'esigne le g\'en\'erateur 
correspondant \`a la $i-$\`eme copie de $E_{2n+2}$).
\end{Ch_n}

\subsection{Les propri\'et\'es de la construction}

V\'erifions maintenant les hypoth\`eses de r\'ecurrence pour le cas $n + 1$.

\begin{Pete}
Du diagramme (\ref{diag1}) on d\'eduit un diagramme de fibrations
\begin{equation}\label{diag2}
\begin{CD}
F_{n+1} @>>>          E       @>>>          G \\
@VVV                  @VVV                  @VVV \\
X_{n+1} @>{i_{n+1}}>> P_{n+1} @>{f_{n+1}}>> G_{n+1} \\
@V{p_{n+1}}VV                  @VVV                  @VVV \\
X_n     @>{i_n}>>     P_n     @>{f_n}>>     G_n
\end{CD}
\end{equation}
De la d\'efinition de $f_{n+1}$ on d\'eduit que
$F_{n+1} \cong K_{4n+3}^{\times (k+h)}$,
donc est $(4n+2)-$connexe et donc
$p_{n+1}^\ast: {\rm H}^\ast X_n \to {\rm H}^\ast X_{n+1}$
est un isomorphisme en degr\'e $i$ tel que $i \leq 4n+2$. $\hfill{\Box}$
\end{Pete}

\begin{Pete}
Si on note par $z_i$ l'\'el\'ement $x_i$ pour $1 \leq i \leq k$
et l'\'el\'ement $y_{i-k}$ pour $k+1 \leq i \leq k+h$,
on a $f_{n+1}^\ast (\iota_{2n+2,i}) = \bar{\iota}_{2n+2,i} + z_i^2$.
Donc $f_{n+1}^\ast (\iota_{2n+2,i}^2) = z_i^4$.
Donc l'image de $z_i$ dans ${\rm H}^\ast X_{n+1}$
est de puissance quatri\`eme nulle.

Ces calculs montrent que tous les g\'en\'erateurs de
$\tilde{\rm H}^\ast X_n$ en degr\'e inf\'erieur ou \'egal \`a
$n+1$ provenant de ${\rm H}^\ast P_n$ sont
d'image nilpotente dans $\tilde{\rm H}^\ast X_{n+1}$. De plus,
\`a cause de la connectivit\'e de la fibre $F_{n+1}$
on n'a pas rajout\'e de g\'en\'erateurs en degr\'e
inf\'erieur ou \'egal \`a $4n+2$. $\hfill{\Box}$
\end{Pete}

Il reste \`a v\'erifier la propri\'et\'e $(3)$.

On consid\`ere le noyau du morphisme d'alg\`ebres de Hopf:
$$f_{n+1}^\ast: {\rm H}^\ast G_{n+1} \to {\rm H}^\ast P_{n+1}.$$
Par l'hypoth\`ese de r\'ecurrence pour $n$, on a
$$(f_n^\ast)^{-1}(Sq_0({\rm PH}^\ast P_n)) \subset Sq_0({\rm PH}^\ast G_n),$$
ce qui implique que le noyau (en tant qu'alg\`ebre de Hopf)
de $f_n^\ast$ est engendr\'e par des cup-carr\'es.
Il faut montrer la m\^eme propri\'et\'e pour $f_{n+1}^\ast$.

\begin{Pete}
On a:
$$(f_{n+1}^\ast)^{-1}(Sq_0({\rm PH}^\ast P_{n+1})) \subset Sq_0({\rm PH}^\ast G_{n+1}).$$
\end{Pete}

\begin{proof}
Consid\'erant la fl\`eche sur les primitifs qui peut s'\'ecrire comme:
$$f_{n+1}^\ast: \oplus_1^{k+h} {\rm F}(2n+2) \bigoplus {\rm PH}^\ast G_n \to
\oplus_1^{k+h} {\rm PH}^\ast E_{2n+2} \bigoplus {\rm PH}^\ast P_n,$$
et qu'on \'ecrira matriciellement
$$\begin{pmatrix}
h& 0 \cr
k&f_{n}^\ast\cr
\end{pmatrix}.$$
Dans cette notation $h$ d\'esigne l'application de
modules instables qui envoie $\iota_{2n+2,i}$ vers $\bar
\iota_{2n+2,i} \in {\rm H}^\ast E_{2n+2}$ et $k$ d\'esigne l'application
de modules instables qui envoie $\iota_{2n+2,i}$ vers $z_i^2$.
Soit $u_1 \in \oplus_1^{k+h} {\rm F}(2n+2) = {\rm PH}^\ast
K_{2n+2}^{\times (k+h)}$ tel que
$h(u_1) \in \oplus_1^{k+h} {\rm PH}^\ast E_{2n+2}$
soit de la forme $Sq_0 v_1$. On veut montrer que $u_1 = Sq_0 w_1$. A
cette fin, on observe que en tant qu'une combinaison lin\'eaire
d'\'el\'ements primitifs
$$u_1 = \sum_{1 \leq t \leq k, s} Sq^{I_{t,s}} \iota_{2n+2,t}$$
est envoy\'ee sur
$$h(u_1)=\sum_{1 \leq t \leq k, s} Sq^{I_{t,s}} \bar{\iota}_{2n+2,t}$$
et que par les propri\'et\'es de $\bar{\iota}_{2n+2}$, cette
combinaison est cup-carr\'e d'un \'el\'ement si et seulement si
elle s'annule. Par cons\'equent, pour tous $t, s$, ${\rm ex}(I_{t,s}) =
2n+2$, c'est-\`a-dire,
$$Sq^{I_{t,s}} \iota_{2n+2,t} = Sq_0 (Sq^{I_{t,s}'} \iota_{2n+2,t})$$
avec $I_{t,s} = (i_1, I_{t,s}')$. Donc on a $u_1 = Sq_0 w_1$ pour
$$w_1 = \sum_{1 \leq t \leq k, s} Sq^{I_{t,s}'} \iota_{2n+2,t}.$$

Soit $(u_1, u_2) \in \oplus_1^{k+h} {\rm F}(2n+2) \bigoplus {\rm PH}^\ast G_n$
tel que $f_{n+1}^\ast (u_1, u_2) \in \oplus_1^{k+h} {\rm PH}^\ast E_{2n+2} \bigoplus {\rm PH}^\ast P_n$ soit de la forme $(Sq_0 v_1, Sq_0 v_2)$. On a donc
$$\begin{array}{rcl}
f_{n+1}^\ast (u_1, u_2)& =& (h(u_1), k(u_1)+f_n^\ast(u_2)) \\
& =& (Sq_0 v_1, Sq_0 v_2).
\end{array}$$
Comme $h(u_1) = Sq_0 v_1$, on a $u_1 = Sq_0 w_1$. Donc
$$f_n^\ast(u_2) = k(u_1) + Sq_0 v_2 = Sq_0 (k(w_1) + v_2)$$
et le r\'esultat d\'ecoule de l'hypoth\`ese de r\'ecurrence.
\end{proof}

\subsection{La cohomologie de $X_\infty$}

Ci-dessus, on a construit les espaces $X_n$ et chaque $X_n$ est la fibre homotopique d'une application $f_n: P_n \to G_n$. Par construction, on a les propri\'et\'es suivantes
\begin{itemize}
\item le noyau (en tant qu'alg\`ebre de Hopf) de $f_n^\ast$ est engendr\'e par des cup-carr\'es;

\item $p_n^\ast: {\rm H}^\ast X_{n-1} \to {\rm H}^\ast X_n$ est un
isomorphisme en degr\'e inf\'erieur ou \'egal \`a $4n-2$;

\item pour tout g\'en\'erateur (en tant qu'alg\`ebre), et donc pour tout \'el\'ement, $x$ de $\tilde{\rm H}^\ast P_{n-1}$ en degr\'e inf\'erieur ou \'egal \`a $n$, soit $\bar{x} = i_{n-1}^\ast(x)$ son image dans $\tilde{\rm H}^\ast X_{n-1}$, alors $p_n^\ast(\bar{x}) \in \tilde{\rm H}^\ast X_n$ est nilpotent.
\end{itemize}

Remarquons que les \'el\'ements non-nilpotents dans
${\rm H}^\ast X_n$ peuvent ne pas provenir de ${\rm H}^\ast P_n$.
Soit $z$ un tel \'el\'ement, dans ce cas le lemme 1 nous garantit
que l'une de ses puissances provient de ${\rm H}^\ast P_n$. Cette
puissance peut alors s'exprimer, modulo un terme nilpotent,  comme
un polyn\^ome en les images des g\'en\'erateurs polyn\^omiaux de
${\rm H}^\ast P_n$. La construction garantit alors que pour $k$
assez grand l'image de $z$ dans ${\rm H}^\ast X_{n+k}$ est
nilpotente.

Par cons\'equent, tout \'el\'ement de
$$\tilde{\rm H}^\ast X_\infty = {\rm colim}_n \tilde{\rm H}^\ast X_n$$
est nilpotent, i.e., $\tilde{\rm H}^\ast X_\infty$ est nilpotente.

\section{Les $2-$polyGEMs stables}

Dans cette section, on donne une d\'emonstration de la conjecture
\ref{GrodalConj} pour les $2-$polyGEMs stables simplement
connexes. Soit $X$ un $2-$polyGEM stable qui est la fibre
homotopique d'un morphisme de H-espaces $f: E \to F$ dont $E$, $F$
sont des $1-$polyGEMs. On va proc\'eder la d\'emonstration en
plusieurs \'etapes selon $E$ et $F$. Rappelons d'abord que
d'apr\`es les r\'esultats de Smith, on a
\begin{equation}\label{smith}
{\rm H}^\ast X \cong {\rm H}^\ast E // {\rm im}(f^\ast) \otimes
{\rm Tor}^*_{{\rm sub-ker}(f^\ast)} ({\Bbb F}_2, {\Bbb F}_2).
\end{equation}

\begin{Cas}
On consid\`ere d'abord le cas o\`u $E$, $F$ sont des produits
d'espaces d'Eilenberg-MacLane ${\rm K}({\Bbb F}_2,n)$. Soit donc $F =
\prod_\alpha {\rm K}({\Bbb F}_2, n_\alpha)$ et $E = \prod_\beta
{\rm K}({\Bbb F}_2, m_\beta)$. Soit $M = \oplus_\alpha {\rm
F}(n_\alpha)$ (resp. $L = \oplus_\beta {\rm F}(m_\beta)$)
l'ensemble des \'el\'ements primitifs de ${\rm H}^\ast F$ (resp.
${\rm H}^\ast E$), le morphisme $f^\ast: {\rm H}^\ast F \to {\rm
H}^\ast E$ induit une application (notée abusivement) $f^*: M \to
L$. Si on note $Q(M) = M/\bar{{\mathcal A}_2}M$ et $Q(L) =
L/\bar{{\mathcal A}_2}L$, on a une application induite
$$Q(f^*): Q(M) \to Q(L).$$
\end{Cas}

{\bf Cas 1.a.} Si $Q(f^*)$ est injective, on sait que $M$ est un
facteur direct dans $L$. Donc $L/M$ est une somme directe de ${\rm
F}(k)$ et
$${\rm H}^\ast E // {\rm im}(f^\ast) \cong {\mathcal U}(L) // {\mathcal U}(M) \cong {\mathcal U}(L/M)$$
contient des \'el\'ements non-nilpotents d\`es que $L \neq M$, ce qui affirme la conjecture \ref{GrodalConj} dans ce cas d'apr\`es la formule (\ref{smith}).

{\bf Cas 1.b.} Si $Q(f^\ast)$ n'est pas injective, alors
on peut trouver (quitte \`a changer la base de $Q(M)$)
un \'el\'ement $\iota_n \in M$ tel que
$$f^\ast(\iota_n) = \sum_{|I|>0} Sq^I \iota_{n-|I|},$$
o\`u ${\rm ex}(I) \leq n-|I|$.

Notons $Q_i$ les d\'erivations de Milnor.
Rappelons que {\bf l'op\'eration $\bf Sq_1$} est d\'efinie sur un module
instable $M$ par $Sq_1 x = Sq^{|x|-1} x$ pour tout $x \in M$.

\begin{Lem}
1. Pour tout $y \in {\mathcal A}_2(Q_0 Q_1 \cdots Q_{n-1} \iota_n)$
on a $f^\ast(y) = 0$.

2. Soit $\omega = Q_0 Q_1 \cdots Q_{n-2} \iota_n$, on a $f^\ast(\omega) = 0$
et si $n \geq 2$, $(Sq_1)^r \omega \neq 0$ pour tout $r \geq 0$.
\end{Lem}

Le sous-module ${\mathcal A}_2(Q_0 Q_1 \cdots Q_{n-1} \iota_n)$ est
isomorphe \`a $\bigwedge^n({\rm F}(1))$, o\`u ${\rm F}(1) \cong
{\mathcal A}_2 u \subset {\Bbb F}_2 [u] = {\rm H}^\ast {\rm B}
{\Bbb F}_2$ (voir \cite{Steenrod}).
Le lemme implique que $\omega \in {\rm sub}$-${\rm ker}(f^\ast)$,
puisque pour tout morphisme $\varphi$ d'alg\`ebres de Hopf r\'eduites primitivement engendr\'ees, $\rm sub$-${\rm ker}(\varphi)$ est l'alg\`ebre de Hopf engendr\'ee par les \'el\'ements primitifs du noyau de $\varphi$.

Gr\^ace au complexe de Koszul, on sait que pour que
${\rm Tor}_{{\rm sub-ker}(f^\ast)}^\ast ({\Bbb F}_2, {\Bbb F}_2)$
contienne des \'el\'ements non-nilpotents, il suffit d'avoir
un \'el\'ement dans $\rm sub$-${\rm ker}(f^\ast)$ sur lequel
l'op\'eration $(Sq_1)^r$ ne s'annule pas pour tout $r \geq 0$.
Comme $\omega \in {\rm sub}$-${\rm ker}(f^\ast)$,
le lemme montre donc que ${\rm Tor}_{{\rm sub-ker}(f^\ast)}^\ast
({\Bbb F}_2, {\Bbb F}_2)$, donc ${\rm H}^\ast X$
contient des \'el\'ements non-nilpotents.

\begin{flushleft}
{\it D\'emonstration du lemme 2.}
En effet, identifiant la classe $\iota_t$ ($t = n -$ $|I|$) avec
$x_1 \cdots x_t \in {\Bbb F}_2[x_1, \cdots, x_t]$ (voir \cite{Steenrod}), 
comme $n-1 \geq t$ le produit de $n$ d\'erivations
de Milnor s'annule sur ${\Bbb F}_2[x_1, \cdots, x_t]$ et on a
\end{flushleft}
$$\begin{array}{rl}
& f^\ast (Q_0 Q_1 \cdots Q_{n-1} \iota_n) \\
=& \sum_{|I|>0} Q_0 Q_1 \cdots Q_{n-1} Sq^I \iota_{n-|I|} \\
=& 0.
\end{array}$$
De plus,
$$Sq_0 (Q_0 Q_1 \cdots Q_{n-2} \iota_n) = Q_0 Q_1 \cdots Q_{n-1} \iota_n$$
et l'application $Sq_0$ est injective dans $L$. Il en r\'esulte que
$$f^\ast (Q_0 Q_1 \cdots Q_{n-2} \iota_n) = 0.$$
Comme $\omega = Q_0 Q_1 \cdots Q_{n-2} \iota_n$ s'identifie \`a
$u \wedge u^2 \wedge \cdots \wedge u^{2^{n-1}}
\in \bigwedge^n({\rm F}(1))$,
on a pour tout $r \geq 0$, $(Sq_1)^r \omega \neq 0$. $\hfill{\Box}$

\begin{Cas}
On consid\`ere le cas g\'en\'eral o\`u $E$, $F$ sont des produits
d'espaces d'Eilenberg-MacLane ${\rm K}({\Bbb F}_2, n)$, ${\rm
K}({\Bbb Z}, n)$ et ${\rm K}({\Bbb Z}/2^h, n)$ ($h, n > 1$). On
remarque d'abord que si on note ${\rm F}'(n+1) = {\mathcal A}_2
(Sq^1 \iota_n)$ le sous-${\mathcal A}_2$-module de ${\rm F}(n)$
engendr\'e par $Sq^1 \iota_n$, on a le

\begin{Lem}
1. ${\rm H}^\ast ({\rm K}({\Bbb Z}, n)) = {\mathcal U}({\rm F}'(n))$
et ${\rm H}^\ast ({\rm K}({\Bbb Z}/2^h, n)) = {\mathcal U}({\rm F}'(n))
\otimes {\mathcal U}({\rm F}'(n+1))$, $h > 1$.

2. ${\rm F}'(n)$ est le module instable librement engendr\'e
par une classe $\iota_n'$ de degr\'e $n$ tel que $Sq^1 \iota_n' = 0$.
\end{Lem}

\begin{proof}
La premi\`ere partie est due \`a J.-P. Serre \cite{Serre}. Alors la deuxi\`eme partie est une cons\'equence du fait que l'ensemble des op\'erations qui s'annulent sur $Sq^1$ est ${\mathcal A}_2 Sq^1$.
\end{proof}

Soit $M = \oplus_\alpha {\rm F}(n_\alpha) \bigoplus \oplus_\gamma {\rm F}'(n_\gamma)$ (resp. $L = \oplus_\beta {\rm F}(m_\beta) \bigoplus \oplus_\delta {\rm F}'(m_\delta)$) l'ensemble des \'el\'ements primitifs de ${\rm H}^\ast F$ (resp. ${\rm H}^\ast E$). Le morphisme $f^\ast: {\rm H}^\ast F \to {\rm H}^\ast E$ induit une application (note abusivement) $f^\ast: M \to L$, et donc une application
$$Q(f^\ast): Q(M) \to Q(L).$$
\end{Cas}

{\bf Cas 2.a.} Si $Q(f^\ast)$ est injective, la conjecture \ref{GrodalConj} est vraie d\`es que $L \neq M$. La d\'emonstration se passe comme dans le cas 1.a, on notera cependant qu'on peut avoir une situation (quitte \`a changer la base) du type
$$Q(f^\ast)(\iota_n') = \iota_n,$$
o\`u le conoyau est ${\rm F}'(n+1)$ qui cr\'ee des \'el\'ements non-nilpotents.

{\bf Cas 2.b.} Si $Q(f^\ast)$ n'est pas injective, on peut supposer
$Q(f^*)$ restreinte \`a $Q(\oplus_\alpha {\rm F}(n_\alpha))$
injective, sinon on est ramen\'e au cas trait\'e plus haut. Alors
on peut se ramener \`a une situation o\`u $M = \oplus_\gamma
{\rm F}'(n_\gamma)$ et trouver (quitte \`a changer la base de $Q(M)$)
un \'el\'ement $\iota_n' \in M$ tel que $Q(f^\ast)(\iota_n') = 0$.
Il est clair que $f^\ast(\iota_n')$ ne contient aucun terme $\iota_n$
car $Sq^1 \iota_n' = 0$ et $Sq^1 \iota_n \neq 0$.

Deux cas peuvent se pr\'esenter, soit $f^\ast(\iota_n')$ contient un
terme $Sq^1 \iota_{n-1}$ et dans ce cas le conoyau de $f^*: M
\to L$ aura un quotient isomorphe \`a ${\rm F}'(n-1)$. On trouve
alors un \'el\'ement non-nilpotent dans le conoyau d\`es que $n \geq 3$.
Si $f^\ast (\iota_n')$ ne contient pas un tel terme, un raisonnement
analogue \`a celui donn\'e plus haut permet de conclure, et de trouver
un \'el\'ement non-nilpotent dans
${\rm Tor}_{{\rm sub-ker}(f^\ast)}^\ast ({\Bbb F}_2, {\Bbb F}_2)$.

\section*{Remerciement}

L'auteur tient \`a remercier Lionel Schwartz pour des discussions
tr\`es importantes sur la formulation de cet article et aussi pour
sa direction.


\begin{thebibliography}{99}
  \bibitem{BP}
    E.H. Brown and F.P. Peterson, {\em Whitehead products and cohomology operations}, Quart. J. Math. Oxford Ser.(2), 15 (1964), pp.116-120.
  \bibitem{Cohen}
    F. Cohen, Communication priv\'ee, 2003.
  \bibitem{Farjoun}
    E. Dror Farjoun, {\em Cellular spaces, null spaces and homotopy localization}, Lecture Notes in Mathematics, 1622. Springer-Verlag, Berlin, (1996).
  \bibitem{FHLT}
    Y. F\'elix, S. Halperin, J.-M. Lemaire and J.-C. Thomas, {\em Mod $p$ loop space homology}, Invent. Math. 95 (1989), no. 2, pp.247-262.
  \bibitem{Grodal}
    J. Grodal, {\em The transcendence degree of the mod $p$ cohomology of finite Postnikov systems}, Stable and unstable homotopy, Toronto, Fields Inst. (1996), pp.111-130.
  \bibitem{Kristensen}
    L. Kristensen, {\em On secondary cohomology operations}, Math. Scand. 12 (1963), pp.57-82.
  \bibitem{LS86In}
    J. Lannes et L. Schwartz, {\em \`A propos de conjectures de Serre et Sullivan}, Invent. Math. 83 (1986), no. 3, pp.593-603.
  \bibitem{LS89Is}
    J. Lannes et L. Schwartz, {\em Sur les groupes d'homotopie des espaces dont la cohomologie modulo $2$ est nilpotente}, Israel J. Math. 66 (1989), no. 1-3, pp.260-273.
  \bibitem{McN}
    C.A. McGibbon and J.A. Neisendorfer, {\em On the homotopy groups of a finite-dimensional space}, Comment. Math. Helv. 59 (1984), no. 2, pp.253-257.
  \bibitem{Milgram}
    J. Milgram, {\em The structure over the Steenrod algebra of some 2-stage Postnikov systems}, Quart. J. Math. Oxford Ser.(2), 20 (1969), pp.161-169.
  \bibitem{MM}
    J.W. Milnor and J.C. Moore, {\em On the structure of Hopf algebras}, Annals of Mathematics(2), 81 (1965), pp.211-264.
  \bibitem{Serre}
    J.-P. Serre, {\em Cohomologie modulo $2$ des complexes d'Eilenberg-MacLane}, Comment. Math. Helv. 27 (1953) pp.198-232.
  \bibitem{Smith}
    L. Smith, {\em The cohomology of stable two stage Postnikov systems}, Illinois J. Math. 11 (1967) pp.310-329.
  \bibitem{Steenrod}
    N.E. Steenrod, {\em Cohomology operations}, Lectures by N.E. Steenrod written and revised by D.B.A. Epstein, Annals of Mathematics Studies, no. 50, Princeton University Press, Princeton, (1962).
\end{thebibliography}
\end{document}